\documentclass[12pt]{article}
\usepackage{amssymb}

\usepackage{color}
\usepackage{colortbl}
\definecolor{lightgray}{gray}{.5}
\definecolor{BR}{rgb}{1,0,0}
\newcommand{\G}{\cellcolor{lightgray}}

\begin{document}
\title{On Two Bijections from $S_n(321)$ to $S_n(132)$}
\author{ Dan Saracino\\Colgate University}
\date{}
\maketitle
\begin{abstract}
In [4], Elizalde and Pak gave a bijection $\Theta:S_n(321)\rightarrow S_n(132)$ that commutes with the operation of taking inverses and preserves the numbers of fixed points and excedances for every $\sigma\in S_n(321).$  In [1] it was shown that another bijection $\Gamma:S_n(321)\rightarrow S_n(132)$ introduced by Robertson in [7] has these same properties, and in [2] a pictorial reformulation of $\Gamma$ was given to make it clearer why $\Gamma$ has these properties. Our purpose here is to give a similar pictorial reformulation of $\Theta$, from which it follows that, although the original definitions of $\Theta$ and $\Gamma$ make them appear quite different, these two bijections are in fact related to each other in a very simple way, by using inversion,  reversal, and complementation.
\end{abstract}
\noindent \textbf{1. Introduction}
\vspace{.25in}

If we write the permutation $\sigma \in S_n$ in one-line notation as $\sigma_1\sigma_2\ldots \sigma_n,$ then a triple of entries $\sigma_i\sigma_j\sigma_k$ with $i< j< k$ is called a \emph{321-pattern} if $\sigma_i> \sigma_j>\sigma_k$ and a \emph{132-pattern} if $\sigma_j> \sigma_k> \sigma_i.$  We denote by $S_n(321)$ (respectively,  $S_n(132)$) the set of all $\sigma \in S_n$ that contain no 321-patterns (respectively, no 132-patterns).

Knuth gave a bijection from $S_n(321)$ to $S_n(132)$ in [5], by giving  bijections from each set onto the set $D_n$ of Dyck paths of length $2n$.  (Knuth's bijection from $S_n(132)$ to $D_n$ is the same as one given later by Krattenthaler [6], although Krattenthaler's definition of it was very different from Knuth's.) A number of other bijections from $S_n(321)$ to $S_n(132)$ have been given since Knuth's (see [3]). In [4], Elizalde and Pak exhibited a bijection $\Theta$ that commutes with the operation of taking inverses and preserves the numbers of fixed points and excedances for every $\sigma \in S_n(321).$ The map $\Theta$ was defined by composing Knuth's bijection from $S_n(321)$ to $D_n$ with the inverse of a modification of Krattenthaler's bijection.

On the other hand, Robertson [7] gave a direct bijection $\Gamma:S_n(321) \rightarrow S_n(132)$ that iteratively replaces the smallest 132-pattern $\sigma_i\sigma_j\sigma_k$ by $\sigma_j\sigma_k\sigma_i$ (``smallest" meaning smallest in the lexicographic ordering of the triples $(i,j,k)$ for all the 132-patterns present). Robertson conjectured at the 2005 Integers conference that $\Gamma$ preserves the number of fixed points for every $\sigma \in S_n(321)$, and Bloom and Saracino proved in [1] that it  preserves the number of fixed points and the number of excedances. The proof depended on first showing, via a somewhat involved argument, that $\Gamma$ also commutes with the operation of taking inverses. Motivated by a desire to find simpler arguments, Bloom and Saracino in [2] gave a pictorial reformulation of $\Gamma$ that freed the results on fixed points and excedances from their dependence on the result about inverses, while rendering all the results more transparent.

Our purpose in the present paper is to give a pictorial reformulation of $\Theta$, from which it follows in particular that, despite the disparity in their original definitions, $\Theta$ and $\Gamma$ are in fact related to each other in a very simple way (using inversion, reversal, and complementation).
\vspace{.25in}

\noindent \textbf{2. Background Definitions and Facts}
\vspace{.25in}

Consider an $n\times n$ array of squares, and represent the square in the $i$th row from the top and the $j$th column from the left by $(i,j).$  Let $T$ be a set of squares in the array, and think of the squares in $T$ as shaded and those not in $T$ as unshaded.
\vspace{.15in}

\noindent\textbf{Definition} [2]. If $\sigma\in S_n$, we say that $T$ is a \emph{template} for $\sigma$ if when (inductively) on each row of the array we place a dot in the first unshaded square from the left that has no dot above it, then the dots are in exactly the squares $(i,\sigma(i)).$
\vspace{.15in}

The pictorial reformulation of $\Gamma$ in [2] was obtained by giving a simple way to turn a certain template for $\sigma\in S_n(321)$ into a template for $\Gamma(\sigma)\in S_n(132).$ The template for $\sigma\in S_n(321)$ was obtained by using what were called the \emph{L-corners} of $\sigma.$
\vspace{.15in}

\noindent\textbf{Definition} [2]. If $\sigma\in S_n(321)$ and $\sigma_i\sigma_j$ is a 21-pattern in $\sigma$ (i.e., if $i< j$ and $\sigma_i> \sigma_j$), then we call $\sigma_i$ a \emph{2-element} and $\sigma_j$ a \emph{1-element}.
\vspace{.15in}

It is clear that if $\sigma\in S_n(321)$ then no entry in $\sigma$ can be both a 2-element and a 1-element. It is also easy to see that the 2-elements in $\sigma$ are in increasing order from left to right, and so are the 1-elements.
\vspace{.15in}

\noindent\textbf{Fact 1} ([2], Lemma 2). Suppose $\sigma\in S_n(321)$ and there are 21-patterns in $\sigma.$ Then
\begin{itemize}
\item[(a)] If $\sigma_i$ is the largest 2-element in $\sigma$ and $\sigma_j$ is the largest 1-element, then 
 $\sigma_i\sigma_j$ is a 21-pattern in $\sigma.$
 \item[(b)] Let $\sigma_s\sigma_t$ be a 21-pattern in $\sigma$ and suppose there exist 21-patterns $\sigma_x\sigma_y$ in $\sigma$ such that $\sigma_x< \sigma_s$ and $\sigma_y< \sigma_t.$ Let $\sigma_i$ and $\sigma_j$ be, respectively, the largest 2-element and the largest 1-element occurring in these patterns $\sigma_x\sigma_y.$ Then $\sigma_i\sigma_j$ is a 21-pattern in $\sigma$.
\end{itemize}
\vspace{.15in}
\noindent\textbf{Definition} [2]. Suppose $\sigma\in S_n(321)$ and there are 21-patterns in $\sigma$. The \emph{first L-corner} of $\sigma$ is the pair $(i,\sigma_j)$, where $\sigma_i$ and $\sigma_j$ are, respectively, the largest 2-element and the largest 1-element in $\sigma.$ If $(s,\sigma_t)$ is the $k$th $L$-corner of $\sigma$ and there are 21-patterns $\sigma_x\sigma_y$ in $\sigma$ such that $\sigma_x< \sigma_s$ and $\sigma_y< \sigma_t$ then the \emph{(k+1)th L-corner} of $\sigma$ is the pair $(i,\sigma_j)$, where $\sigma_i$ and $\sigma_j$ are, respectively, the largest 2-element and the largest 1-element occurring in these patterms $\sigma_x\sigma_y.$
\vspace{.15in}

 If the $L$-corners of $\sigma$ are $(p_1,v_1),\ldots, (p_t,v_t)$ (``$p$" for ``position", ``$v$" for ``value"), then for each $(p_i,v_i)$ let $T_i$ be the set of squares $$T_i=\{(p_i,j):j\leq v_i\}\cup \{(j,v_i): j\leq p_i\}.$$ The squares in $T_i$ form a reversed $L$ having its corner at $(p_i,v_i)$ and extending horizontally to the left-hand border of the array and vertically to its top border. Let $T_{\sigma}$ be the union of the $T_i$'s, so that $T_{\sigma}$ is the union of a set of nested reversed $L$'s. (Note that the smallest reversed $L$ may degenerate to a line segment.)
 \vspace{.15in}
 
 \noindent\textbf{Fact 2} ([2], Lemma 3). For every $\sigma\in S_n(321)$, the set $T_{\sigma}$ is a template for $\sigma.$
 \vspace{.15in}
 
 If the $L$-corners of $\sigma$ are $(p_1,v_1),\ldots,(p_t,v_t)$ with $p_t<\ldots < p_1$ and $v_t<\ldots <\ v_1$ then let $\widehat{T_{\sigma}}$ be the set of squares obtained from $T_{\sigma}$ by replacing each $T_i$ by an inverted $L$ having its corner at $(i,i)$ and vertical and horizontal legs consisting of $p_i$ and $v_i$ squares, respectively (including the corner in each leg).
 \vspace{.15in}
 
 \noindent\textbf{Fact 3} ([2], Theorem 1).  For every $\sigma\in S_n(321),$ $\widehat{T_{\sigma}}$ is a template for $\Gamma(\sigma).$
 \vspace{.15in}
 
 \noindent\textbf{Example}. If $\sigma\in S_8(321)$ is 14237586 then the $L$-corners of $\sigma$ are (7,6), (5,5), and (2,3), so $T_{\sigma}$ is

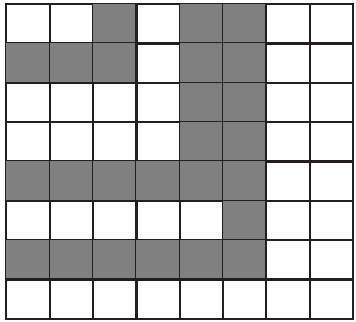
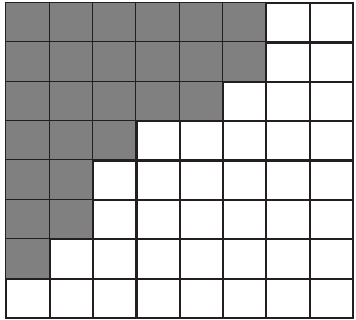
\begin{center}
\begin{tabular}{|p{4pt}|p{4pt}|p{4pt}|p{4pt}|p{4pt}|p{4pt}|p{4pt}|p{4pt}|}
\hline
&&\G&&\G&\G&&\\
\hline
\G&\G&\G&&\G&\G&&\\
\hline
&&&&\G&\G&&\\
\hline
&&&&\G&\G&&\\
\hline
\G&\G&\G&\G&\G&\G&&\\
\hline
&&&&&\G&&\\
\hline
\G&\G&\G&\G&\G&\G&&\\
\hline
&&&&&&&\\
\hline
\end{tabular}

\end{center}
and $\widehat{T_{\sigma}}$ is

\begin{center}
\begin{tabular}{|p{4pt}|p{4pt}|p{4pt}|p{4pt}|p{4pt}|p{4pt}|p{4pt}|p{4pt}|}
\hline
\G&\G&\G&\G&\G&\G&&\\
\hline
\G&\G&\G&\G&\G&\G&&\\
\hline
\G&\G&\G&\G&\G&&&\\
\hline
\G&\G&\G&&&&&\\
\hline
\G&\G&&&&&&\\
\hline
\G&\G&&&&&&\\
\hline
\G&&&&&&&\\
\hline
&&&&&&&\\
\hline
\end{tabular}
\end{center}
so $\Gamma(\sigma)$ is 78643521.
\vspace{.15in}

Since we want to obtain a similar pictorial reformulation of $\Theta$, we recall the original definition from [4].

For $\sigma\in S_n(321),$ we first use the Robinson-Schensted-Knuth correspondence (see, for example, [8]) to obtain the insertion tableau and recording tableau for $\sigma.$ Each tableau will consist of the elements of $\{1,2,\ldots,n\}$ arranged in two rows (unless $\sigma $ is the identity element of $S_n$, in which case each will have one row).  The insertion tableau is formed in stages, by using the entries of $\sigma$ from left to right.  Assuming that $\sigma_1,\ldots,\sigma_{i-1}$ have been placed in the insertion tableau, we place $\sigma_i$ as follows. If $\sigma_i$ is larger than all the elements that are in the first row, then we place $\sigma_i$ at the right end of the first row. If, on the other hand, $\sigma_j$ is the leftmost element in the first row that is larger than $\sigma_i$, we replace $\sigma_j$ by $\sigma_i$ and place $\sigma_j$ at the right end of the second row. (We then say that $\sigma_i$ has ``bumped" $\sigma_j$ to the second row.)  The recording tableau for $\sigma$ has the same shape as the insertion tableau and is obtained by placing each $i$ in the position that became occupied for the first time when $\sigma_i$ was placed in the insertion tableau.
\vspace{.15in}

\noindent\textbf{Example}. If $\sigma$ is 14237586, the insertion and recording tableaux are\\

\vspace{.15in}

\(\begin{array}{ccccc}
1&2&3&5&6\\
4&7&8
\end{array}\)
\vspace{.15in}

\noindent and
\vspace{.15in}

\(\begin{array}{ccccc}
1&2&4&5&7\\
3&6&8
\end{array}\).
\vspace{.15in}

To continue with the definition of $\Theta(\sigma)$, we next form a Dyck path of length $2n$ by getting the first $n$ steps from the insertion tableau and the last $n$ steps from the recording tableau. For the first half, for each $i\in \{1,\ldots,n\}$ we take an upstep ($u$) if $i$ is in the first row of the insertion tableau and a downstep ($d$) if $i$ is in the second row.  For the second half, we form a sequence of $u$'s and $d$'s in the same way from the recording tableau, and then write this sequence in reverse order and interchange $u$'s and $d$'s. We append the result to the first half to complete the Dyck path of length $2n.$
\vspace{.15in}

\noindent\textbf{Example}.  If $\sigma$ is 14237586, the first half of the Dyck path is $uuuduudd$ and the second half is $ududdudd$.
\vspace{.15in}

To obtain $\Theta(\sigma)$, we start with an $n\times n$ array of squares.  We draw a path consisting of edges of the squares, by starting at the point in the lower left-hand corner of the array and (reading the Dyck path from left to right) going up one edge for each $u$ and right one edge for each $d$.
\vspace{.15in}

\noindent\textbf{Example}.  If $\sigma$ is 14237586 then we obtain the path  that is the border between the shaded and unshaded regions in

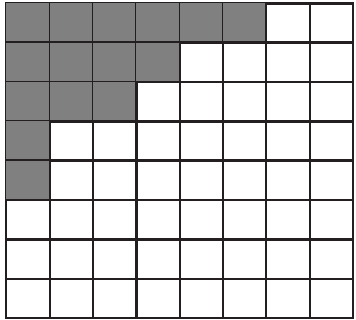
\begin{center}
\begin{tabular}{|p{4pt}|p{4pt}|p{4pt}|p{4pt}|p{4pt}|p{4pt}|p{4pt}|p{4pt}|}
\hline
\G&\G&\G&\G&\G&\G&&\\
\hline
\G&\G&\G&\G&&&&\\
\hline
\G&\G&\G&&&&&\\
\hline
\G&&&&&&&\\
\hline
\G&&&&&&&\\
\hline
&&&&&&&\\
\hline
&&&&&&&\\
\hline
&&&&&&&\\
\hline
\end{tabular}

\end{center}

\vspace{.15in}
Finally, $\Theta(\sigma)$ is defined to be the element of $S_n(132)$ whose diagram is the portion of the array that is to the left of the path.  As always for elements of $S_n(132),$  the diagram is a template for the permutation.

It will be helpful in Section 3 to note that the second half of the path drawn in the last step of the definition of $\Theta(\sigma)$ can be obtained from the recording tableau for $\sigma$ by starting at the point in the upper right-hand corner of the array and, for each $i\in \{1,\ldots,n\}$, going left one edge if $i$ is in the first line of the tableau and down one edge if $i$ is in the second line.
\vspace{.25in}

\noindent\textbf{3. $\Theta$ and its Relationship to $\Gamma$}

\vspace{.25in}
\noindent\textbf{Notation}.  Fix $n$.  For any $i\in \{1,\ldots,n\},$ let $\overline{i}=n+1-i.$
\vspace{.15in}

\noindent{\textbf{Lemma 1}. Let $\sigma\in S_n(321)$ and suppose the second rows of the insertion and recording tableaux for $\sigma$ are $$
a_1\  a_2\ \ldots \ a_n \ \ \textrm{and}\ \ b_1\  b_2\  \ldots\  b_n$$
respectively.  Then the template for $\Theta(\sigma)$ produced by the definition of $\Theta(\sigma)$ consists of $k$ inverted $L$'s, $L_1,\ldots, L_k$, such that the corner of $L_i$ is the square $(i,i)$ and the vertical and horizontal legs of $L_i$ consist of $\overline{a_i}$ and $\overline{b_i}$ squares, respectively (including the corner in each leg).
\vspace{.15in}

\noindent\emph{Proof}. Considering the construction of the path of edges of squares in the definition of $\Theta(\sigma),$ we see that the horizontal edges in the first half of the path occur in the leftmost $k$ columns of the array, and that for $1 \leq i \leq k$  the $i$th horizontal edge is $a_i-i \leq n-i$  edge-lengths above the bottom border of the array. Likewise, viewing the second half of the path as originating at the upper right-hand corner of the array, the vertical edges in the second half occur in the top $k$ rows of the array, and for $1\leq i\leq k$ the $i$th vertical edge from the top is $b_i-i \leq n-i$ edge-lengths to the left of the right border of the array. Therefore the template consisting of those squares to the left of the path is comprised of $k$ inverted $L$'s, $L_1,\ldots,L_k$, such that the corner of $L_i$ is the square $(i,i).$

Now consider the vertical leg of $L_i$. This leg consists of those squares in the lower $n-i+1$ rows of the array that are above the $i$th horizontal edge in the first half of the path. Since this edge is $a_i-i$ edge-lengths above the bottom border of the array, the number of squares in the vertical leg of $L_i$ (including the corner) is $$(n-i+1)-(a_i-i)=\overline{a_i}.$$ Viewing the second half of the path as originating at the upper right-hand corner of the array, we see in the same way that the horizontal leg of $L_i$ consists of $\overline{b_i}$ squares.   $\Box$
\vspace{.15in}

We next need an easy way to read the second rows of the tableaux for $\sigma$ directly from $\sigma.$
\vspace{.15in}

\noindent\textbf{Notation}.  If $\sigma$ is the permutation $\sigma_1\sigma_2\ldots \sigma_n$, then $r\sigma$ denotes the permutation $\sigma_n\sigma_{n-1}\ldots \sigma_1$ and $c\sigma$ denotes the permutation $\overline{\sigma_1}\hspace{.05in}\overline{\sigma_2}\ldots \overline{\sigma_n}.$  We denote the inverse of $\sigma$ by $i\sigma.$
\vspace{.15in}

\noindent\textbf{Definition}.  If $\sigma\in S_n(321)$ then the \emph{rcL-corners} of $\sigma$ are $(v_1,p_1),\ldots,(v_t,p_t)$, where $(\overline{p_1},\overline{v_1}),\ldots, (\overline{p_t},\overline{v_t})$ are the $L$-corners of $rc\sigma.$ 
\vspace{.15in}

In this definition, each $v_i$ is a 2-element and each $p_i$ is the position of a 1-element.
\vspace{.15in}

\noindent\textbf{Example}.  If $\sigma$ is 14237586 then $rc\sigma$ is 31426758, the $L$-corners of which are $(6,5),(3,2),$ and $(1,1),$ so the $rcL$-corners of $\sigma$ are $(4,3),(7,6),$ and $(8,8).$
\vspace{.15in}

\noindent\textbf{Lemma 2.} If $\sigma\in S_n(321)$ and the $rcL$-corners of $\sigma$ are $(v_1,p_1),\ldots, (v_t,p_t)$ in increasing order, then
\begin{itemize}
\item[(a)] The smallest 2-element in $\sigma$ is $v_1$, the smallest 1-element is $\sigma_{p_1}$, and $v_1\sigma_{p_1}$ is a 21-pattern
\item[(b)] If $1\leq j\leq t$ and $S$ is the set of all 21-patterns $\sigma_x\sigma_y$ in $\sigma$ such that $\sigma_x> v_j$ and $\sigma_y> \sigma_{p_j}$, then $v_{j+1}$ is the smallest 2-element occurring in the members of $S$ and $\sigma_{p_{j+1}}$ is the smallest 1-element, and $v_{j+1}\sigma_{p_{j+1}}$ is a 21-pattern.  There are no 21-patterns $\sigma_x\sigma_y$ in $\sigma$ such that $\sigma_x> v_t$ and $\sigma_y> \sigma_{p_t}.$
\end{itemize}
\vspace{.15in}

\noindent\emph{Proof}.  This follows from Fact 1 and the definition of the $rcL$-corners.  $\Box$
\vspace{.15in}

\noindent\textbf{Lemma 3}. If $\sigma\in S_n(321)$ and the $rcL$-corners of $\sigma$ are $(v_1,p_1),\ldots, (v_t,p_t)$ in increasing order, then the second rows of the insertion and recording tableaux for $\sigma$ are $$v_1\ v_2\ \ldots\ v_t\ \ \textrm{and}\ \ p_1\ p_2\ \ldots\ p_t,$$ respectively.
\vspace{.15in}

\noindent\emph{Proof}. We must show that $\sigma_{p_i}$ ``bumps" $v_i$ for $1\leq i \leq t$, and that these are the only ``bumps" that occur in the construction of the insertion and recording tableaux. Note that if $\sigma_s$ bumps $\sigma_t$ then $\sigma_t$ is a 2-element and $\sigma_s$ is a 1-element.
\vspace{.1in}

\noindent\emph{Claim}:  Fix $1\leq j \leq t.$ Then $\sigma_{p_i}$ bumps $v_i$ for $1\leq i \leq j.$ If $\sigma_k\leq v_j$ and $\sigma_k$ gets bumped during the construction of the tableaux for $\sigma$, then $\sigma_k\in \{v_1,\ldots, v_j\}.$ If $\sigma_{\ell}\leq \sigma_{p_j}$ and $\sigma_{\ell}$ bumps some $\sigma_z$ during the construction of the tableaux, then $\sigma_{\ell}\in \{\sigma_{p_1},\ldots, \sigma_{p_j}\}.$
\vspace{.15in}

If we can establish this claim then Lemma 3 will be proved, for there cannot exist $\sigma_x> v_t$ and $\sigma_y> \sigma_{p_t}$ such that $\sigma_y$ bumps $\sigma_x.$ (Any such $\sigma_x$ and $\sigma_y$ would yield a 21-pattern $\sigma_x\sigma_y$ that would contradict the last sentence of Lemma 2(b).)

We will prove the claim by induction on $j$.

For $j=1$, note that no $\sigma_k< v_1$ can ever get bumped, because no such $\sigma_k$ is a 2-element, by Lemma 2(a). Similarly, no $\sigma_{\ell}< \sigma_{p_1}$ can bump any $\sigma_z$ because no such $\sigma_{\ell}$ is a 1-element. Since $v_1$ is placed in the first row of the insertion tableau before $\sigma_{p_1}$ is, and no $\sigma_{\ell}< \sigma_{p_1}$ can bump $v_1$, it follows from the fact that the elements of the first row of the tableau are in increasing order from left to right that $\sigma_{p_1}$ bumps $v_1$, because $v_1$ is to the left of all the other 2-elements in the first row.

Now assume that the claim holds for $1\leq j < t.$ Note that if $\sigma_x> v_j$ and $\sigma_y> \sigma_{p_j}$ and $\sigma_y$ bumps $\sigma_x$, then since $\sigma_x\sigma_y$ is a 21-pattern, we have $\sigma_x\geq v_{j+1}$ and $\sigma_y\geq \sigma_{p_{j+1}}$ by Lemma 2(b).  It remains to show that $\sigma_{p_{j+1}}$ bumps $v_{j+1}.$  But by the induction hypothesis and what we have already said, no $\sigma_y< \sigma_{p_{j+1}}$ bumps $v_{j+1}$, and no $\sigma_x< v_{j+1}$ is bumped by $\sigma_{p_{j+1}}.$  $\Box$
\vspace{.15in}

\noindent\textbf{Theorem 1}.  If $\sigma\in S_n(321)$ and the $rcL$-corners of $\sigma$ are $(v_1,p_1),\ldots, (v_t,p_t)$ in increasing order then we get a template for $\Theta(\sigma)$ by taking inverted $L$'s, $L_i,\ldots, L_t$ such that the corner of $L_i$ is the square $(i,i)$ and the vertical and horizontal legs of $L_i$ consist of $\overline{v_i}$ and $\overline{p_i}$ squares, respectively.
\vspace{.15in}

\noindent\emph{Proof}.  This follows from Lemmas 1 and 3.  $\Box$
\vspace{.15in}

To give an analog of Fact 3 for $\Theta$, we need an appropriate pictorial representation for $\sigma.$ We will use a modification of the notion of a template.
\vspace{.15in}

\noindent\textbf{Definition}. If $T$ is a subset of an $n\times n$ array and$$\overline{T}=\{(\overline{i},\overline{j}):\ (i,j)\in T\}$$ then we say that $T$ is an \emph{rc-template} for $\sigma\in S_n$ if $\overline{T}$ is a template for $rc\sigma.$
\vspace{.15in}

So $T$ is an $rc$-template for $\sigma$ if when we start with the bottom row and (inductively) on each row we place a dot in the first unshaded square from the right that has no dot below it, then the dots are in exactly the squares $(i,\sigma(i).)$

It is immediate that if $\sigma\in S_n(321)$ and $T_{rc\sigma}$ is the template for $rc$$\sigma$ specified in Fact 2 then $\overline{T_{rc\sigma}}$ is an $rc$-template for $\sigma.$ If the $rcL$-corners of $\sigma$ are $(v_1,p_1),\ldots, (v_t,p_t)$, then $\overline{T_{rc\sigma}}$ is the union of $t$ inverted $L$'s having their corners at the squares $(p_i,v_i)$ and extending to the right-hand and bottom borders of the array.
\vspace{.15in}

\noindent\textbf{Example}.  If $\sigma$ is 14237586 then $\overline{T_{rc\sigma}}$ is
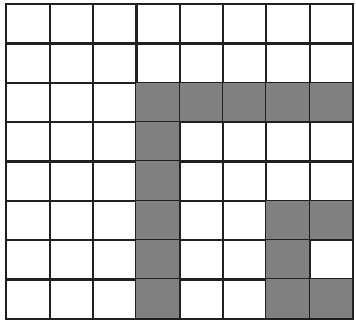
\vspace{.15in}

\begin{center}
\begin{tabular}{|p{4pt}|p{4pt}|p{4pt}|p{4pt}|p{4pt}|p{4pt}|p{4pt}|p{4pt}|}
\hline
&&&&&&&\\
\hline
&&&&&&&\\
\hline
&&&\G&\G&\G&\G&\G\\
\hline
&&&\G&&&&\\
\hline
&&&\G&&&&\\
\hline
&&&\G&&&\G&\G\\
\hline
&&&\G&&&\G&\\
\hline
&&&\G&&&\G&\G\\
\hline
\end{tabular}
\end{center}
\vspace{.15in}

\noindent\textbf{Theorem 2}.  If $\sigma\in S_n(321)$ then we obtain a template for $\Theta(\sigma)$ by taking $\overline{T_{rc\sigma}}$, sliding the largest inverted $L$ so that its vertex becomes the square $(1,1)$, then sliding the second largest so that its vertex becomes the square $(2,2)$, and so on, and then finally flipping all the inverted $L$'s across the diagonal from upper left to lower right.
\vspace{.15in}

\noindent\emph{Proof}. This follows from Theorem 1, since in $\overline{T_{rc\sigma}}$ the inverted $L$ corresponding to $(v_i,p_i)$ has horizontal and vertical legs consisting of $\overline{v_i}$ and $\overline{p_i}$ squares, respectively.  $\Box$
\vspace{.15in}

\noindent\textbf{Theorem 3}.  For every $\sigma\in S_n(321)$ we have $\Theta(\sigma)=\Gamma(irc\sigma).$
\vspace{.15in}

\noindent\emph{Proof}.  If $\sigma$ has $rcL$-corners $(v_1,p_1),\ldots, (v_t,p_t)$  then $T_{rc\sigma}$ is the union of reversed $L$'s, $L_1,\ldots, L_t$, where $L_i$ has its corner at $(\overline{p_i},\overline{v_i})$ and extends to the left-hand and top borders of the array.  $T_{irc\sigma}$ is obtained from $T_{rc\sigma}$ by replacing each corner $(\overline{p_i},\overline{v_i})$ by $(\overline{v_i},\overline{p_i}).$  The result now follows from Theorem 1 and Fact 3.  $\Box$
\vspace{.15in}

By Theorem 3, proving that one of $\Gamma$ or $\Theta$ commutes with the operation $i$, or preserves the number of fixed points or excedances, immediately yields the same result for the other of $\Gamma$ or $\Theta.$
\vspace{.25in}

\centerline{\textbf{References}}
\vspace{.15in}

\noindent1. J. Bloom and D. Saracino, On bijections for pattern-avoiding permutations, \emph{J. Combin. Theory Ser. A} \textbf{116} (2009), 1271-1284.\\
2. J. Bloom and D. Saracino, Another look at pattern-avoiding permutations, \emph{Adv. Appl. Math.} \textbf{45} (2010), 395-409.\\
3. A. Claesson and S. Kitaev, Classification of bijections between 321- and 132-avoiding permutations, \emph{S\'{e}m. Lothar. de Combin.} \textbf{60} (2008), B60d.\\
4. S. Elizalde and I. Pak, Bijections for refined restricted permutations, \emph{J. Combin. Theory Ser. A} \textbf{105} (2004), 207-219.\\
5. D. Knuth, The Art of Computer Programming, Vol. III, Addison-Wesley, Reading, MA, 1973.\\
6. C. Krattenthaler, Permutations with restricted patterns and Dyck paths, \emph{Adv. Appl. Math.} \textbf{27} (2001), 510-530.\\
7.A Robertson, Restricted permutations from Catalan to Fine and back, \emph{S\'{e}m. Lothar. de Combin.} \textbf{50} (2004), B50g.\\
8. R. Stanley, Enumerative Combinatorics, Vol. I, II, Cambridge Univ. Press, Cambridge, 1997, 1999.

\end{document}